\documentclass[preprint,12pt]{elsarticle}

\usepackage{amssymb}
\usepackage{amsmath}
\usepackage[ruled,vlined]{algorithm2e}

\journal{arxiv.org}

\begin{document}

\begin{frontmatter}



\title{Nonlinearity continuation method for steady-state groundwater flow modeling in variably saturated conditions}


\author[label1,label2]{Denis Anuprienko}
\author[label1,label2]{Ivan Kapyrin}
\address[label1]{Marchuk Institute of Numerical Mathematics, Russian Academy of Sciences}
\address[label2]{Nuclear Safety Institute, Russian Academy of Sciences}


\begin{abstract}
Application of nonlinearity continuation method to numerical solution of steady-state groundwater flow in variably saturated conditions is presented. In order to solve the system of nonlinear equations obtained by finite volume discretization of steady-state Richards equation, a series of problems with increasing nonlinearity are solved using the Newton method. This approach is compared to pseudo-transient method on several test cases, including real site problems and involving parallel computations.
\end{abstract}


\end{frontmatter}


\section{Introduction}
\label{}
Numerical modeling of groundwater flow in variably saturated soils is an important part of solution of various hydrogeological problems involving vadose zone. These problems arise in environmental sciences and include safety assessment of radioactive waste repositories and landfills, water resources management and fertilizers application optimization. 
Groundwater flow in variably saturated conditions is governed by nonlinear partial differential Richards equation \cite{richards1931capillary,bear2010modeling} with additional constitutive relationships between water pressure head, water content and relative permeability of the soil. Existing analytical and semi-analytical solutions \cite{parlange1999analytical} of Richards equation are valid for limited cases, and numerical solution is the common choice. Nonlinearity of the equation makes numerical solution a challenging task with numerous approaches having been developed, an overview of which can be found in \cite{farthing2017numerical}.
\par
In addition to transient problems, steady-state problems are also of interest. Discretization of steady-state Richards equation leads to a system of nonlinear equations which then can be solved by an iterative method. Nonlinear solvers such as Newton method can experience severe convergence problems since obtaining a suitable initial guess is not trivial. Another approach is the pseudo-transient method \cite{farthing2003efficient}, which introduces time to the system and lets it evolve until equilibrium is reached, leading to solution of transient Richards equation. For stability reasons discretization in time usually employs fully implicit scheme which again results in a system of nonlinear equations at each time step, although now with ability to change time step size to improve convergence. Nonlinear solvers (modified Picard \cite{celia1990general} and Newton with techniques such as primary variable switching \cite{diersch1999primary}) within the pseudo-transient method can still exhibit poor convergence in case of dry and heterogeneous soils or highly nonlinear constitutive relationships, leading to a large number of time steps and overall high computational complexity. Therefore, it is desirable to have a solution approach which needs no time-stepping and works directly with discretized steady-state Richards equation. As mentioned before, the main problem is the lack of initial guess sufficiently close to solution. Techniques like line search \cite{farthing2003efficient} and combined Picard-Newton approach \cite{paniconi1994comparison} can improve convergence, but still can fail for complex problems.
\par
To overcome difficulties related to lack of suitable initial guess for Newton method a nonlinearity continuation method can be applied. In methods of this type a suitable initial guess for Newton method is obtained through solution of a series of other problems with incrementally varying nonlinearity \cite{continCFD}. This approach for steady-state Richards equation is implemented in GeRa \cite{kapyrin2015integral} hydrogeological modelling software, the authors being part of the GeRa development team. The nonlinearity continuation method is compared with the pseudo-transient method on several modeling cases.
\par
The article is organized as follows. In the second section the mathematical model of steady-state groundwater flow in unconfined conditions is formulated. In the third section the numerical methods used to discretize the problem and solve the arising nonlinear systems are stated, the algorithms of Newton method with line search and nonlinearity continuation method applied to the problem are presented. In the fourth section the proposed numerical technique is tested on a couple of problems featuring heterogeneous and anisotropic media properties with two types of finite volume discretization in space. Conclusions are provided in the end.

\section{Governing equations}
\label{}
Steady-state groundwater flow in variably saturated conditions is described by steady-state Richards equation

\begin{equation}\label{eq:RichStat}
-\nabla \cdot (K_r(h)\mathbb{K}\nabla h) = Q,
\end{equation}

which is derived from Richards equation \cite{richards1931capillary} for transient flow, which can be formulated with additional storage term \cite{bear2010modeling} as:

\begin{equation}\label{eq:Richards}
\frac{\partial \theta(h)}{\partial t} + s_{stor}S(h)\frac{\partial h}{\partial t} -\nabla \cdot (K_r(h)\mathbb{K}\nabla h) = Q.
\end{equation}
Here the following variables are used:
\begin{itemize}
\item $h$ -- hydraulic head, related to capillary pressure head $\Psi$ as $h = \Psi + z$;
\item $\theta(h)$ -- volumetric water content in medium;
\item $S(h)$ -- water saturation in medium;
\item $s_{stor}$ -- specific storage coefficient;
\item $\mathbb{K}$ -- hydraulic conductivity tensor, a 3$\times$3 s.p.d. matrix;
\item $K_r(h)$ -- relative permeability for water in medium;
\item $Q$ -- specific sink and source terms.
\end{itemize}

In this paper the main goal is to solve steady-state equation \eqref{eq:RichStat}, although transient equation \eqref{eq:Richards} is considered within a pseudo-transient method.
\par

Expressions relating relative permeability $K_r(h)$, water content $\theta(h)$ and hydraulic head $h$ are required in addition to equations  \eqref{eq:RichStat} and \eqref{eq:Richards}. In this paper  piecewise linear dependencies proposed in \cite{anuprienko2018modeling} are considered. These functions are designed for numerical solution with the finite volume method and depend on spatial discretization, possibly varying from cell to cell.
Namely, the dependence of water content $\theta_E$ in a cell $E$ on hydraulic head $h_E$ in that cell is defined as follows
\begin{equation}\label{eq:MC}
\theta_E(h_E) = 
\begin{cases}
\phi, h_E > h_{E,\max},\\
\phi \cdot \frac{h_E - h_{E,\min}}{h_{E,\max} - h_{E,\min}}, h_{E,r} < h_E \le h_{E,\max},\\
\phi \cdot (\alpha_{\phi} - \alpha_{\theta}(h_{E,r} - h_E)), h_E \le h_{E,r},
\end{cases}
\end{equation}
where $\phi$ is the medium porosity, $h_{E,\max}$ and $h_{E,\min}$ are the maximal and minimal vertical node coordinates of the cell $E$ 
and $h_{E,r}$ is such that water content calculated by the second linear part in \eqref{eq:MC} is equal to $\phi \cdot \alpha_{\phi}$, namely,
\begin{equation}
h_{E,r} = h_{E,\min} + \alpha_{\phi}(h_{E,\max} - h_{E,\min}).
\end{equation}
Therefore the model needs two parameters: $\alpha_{\phi}$ and $\alpha_{\theta}$ which should be small enough.
\par
Relative permeability for a cell $E$ is assumed to be equal to the saturation:
\begin{equation}\label{eq:Kr}
K_{r,E}(h_E) = S(h_E) = \frac{\theta_E(h_E)}{\phi}.
\end{equation}
\par 
Compared to commonly used functions proposed by van Genuchten \cite{van1980closed} and Mualem \cite{mualem1976new} that are constructed for pressure head rather than hydraulic head, piecewise linear functions require smaller number of media parameters (which are often hard to obtain for real-life problems) and tend to more sharp water table boundary at the cost of ability to predict capillary effects.


\section{Approaches to solution of nonlinear problems}\label{label}
The considered nonlinear system arises from finite volume discretization of Richards equation on unstructured meshes. Two types of unstructured meshes are used: the first consists mainly of triangular prisms with occurence of tetrahedra or pyramids in cases like geological layer pinch-outs; the second includes octree-based hexahedral grids with cut cells on domain boundary \cite{plenkin2015adaptive}. Finite volume schemes used employ different flux approximations across cell faces:
\begin{itemize}
\item conventional linear two-point flux approximation (TPFA);
\item linear multipoint flux approximation (MPFA), O-scheme \cite{aavatsmark1998discretization}.
\end{itemize}

\subsection{Solution of nonlinear systems}
Due to nonlinearity in equation \eqref{eq:RichStat}, its discretization results in a system of nonlinear equations with respect to the grid unknowns of the form

\begin{equation}\label{eq:NonlinSys}
F(h) = 0,
\end{equation}

which is to be solved numerically. For this purpose iterative solvers such as Newton, Picard or their combination \cite{farthing2003efficient}\cite{paniconi1994comparison} can be used. In this work the Newton method is used, so that at $(k+1)$-th iteration a linear system

\begin{equation}\label{eq:JacSys}
\textbf{J}(h^k)\Delta h = -F(h^k),
\end{equation}

is created and solved, where $\textbf{J} = [\partial F_i / \partial h_j]_{ij}$ is the Jacobian matrix and $\Delta h$ is the update vector; then the update procedure is performed:

\begin{equation}
h^{k+1} = h^k + \Delta h.
\end{equation}

After the update, a correction procedure described in \cite{anuprienko2018modeling} is performed that takes into account use of expressions \eqref{eq:MC} and \eqref{eq:Kr} and ensures mass-conservative changes in hydraulic head; it was also noted to improve convergence of nonlinear solvers for transient problems.
\par
 The iterations continue until for the residual $r_k = F(h_k)$ its Euclidean norm is sufficiently reduced compared to the initial residual norm or maximal absolute value of its elements becomes small enough, that is, until one of the following conditions is satisfied:

\begin{equation}
||r_k||_2 < \varepsilon_{rel}\cdot ||r_0||_2,
\end{equation}
\begin{equation}
||r_k||_{\infty} < \varepsilon_{abs}.
\end{equation}

Sometimes it is useful to perform relaxation procedure \cite{paniconi1994comparison}, in which obtained solution is replaced by a convex linear combination of it with the solution from the previous iteration. A \textit{relaxation parameter} $\Omega: 0 < \Omega \le 1$ is chosen and the update procedure is changed to this:

\begin{equation}
h^* = h^k + \Delta h,
\end{equation}
\begin{equation}\label{eq:RelaxUpdate}
h^{k+1} = \Omega h^* + (1-\Omega) h^k = h^k + \Omega \Delta h.
\end{equation}

It can be seen from \eqref{eq:RelaxUpdate} that relaxation is equivalent to rescaling the update vector $\Delta h$ in the update procedure.
\par
Relaxation can prevent the Newton method from oscillating between several values \cite{paniconi1994comparison}, but it demolishes its quadratic convergence. Therefore, it is reasonable to use relaxation during first iterations and at the end converge to the solution quadratically. However, it is not clear how to choose that number of iterations and the relaxation parameter, and we turn to automatic approach of choosing relaxation parameter, a simple version of line search (more sophisticated versions of which can be found in \cite{farthing2003efficient}) that is described in algorithm \ref{alg:NewtonAutoRelax}.

\begin{algorithm}[H]\label{alg:NewtonAutoRelax}
\SetAlgoLined
    Set initial hydraulic head values $h_0$\;
    $k = 0$\;
    \While{$k < \text{maxiter}$}{
        \If{$k = 0$}{
            Compute $r_0$, $||r_0||_2$, $||r_0||_{\infty}$\;
        }

        Create and solve $\textbf{J}(h^k)\Delta h = -F(h^k)$\;
        \eIf{$k \ge 5$}{
            $\Omega = 1$\;
            \While{$\Omega > \Omega_{min}$}{
                Compute $r = F(h^k + \Omega \Delta h)$\;
                \eIf{$||r||_2 < ||r_k||_2$}{
                    $h^{k+1} = h^k + \Omega \Delta h$\;
                    break\;
                }{
                   $\Omega = \gamma \cdot \Omega$\;
                }
            }
            \If{didn't find $\Omega$}{
                failed\;
                stop\;
            }
        }{
            $h^{k+1} = h^k + \Delta h$\;
        }
        Perform correction of $h^{k+1}$ (see \cite{anuprienko2018modeling})\;
        $r_{k+1} = F(h^{k+1})$\;

        \If{$||r_{k+1}||_2 < \varepsilon_{rel} \cdot ||r_0||_2$ or $||r_{k+1}||_{\infty} < \varepsilon_{abs}$}{
            converged\;
            stop\;
        }
    $k = k + 1$\;
    }
 \caption{Newton method with line search}
\end{algorithm}
In this algorithm $\gamma$ is the decreasing factor, which is chosen to be $\gamma = 0.25$. $\Omega_{min}$ is chosen such that 7 iterations of $\Omega$ refinement can be done. Note that line search is not applied for the first 5 iterations. This approach is based on following observations: first, that finding suitable $\Omega$ at the first iterations is often not possible; second, that Newton method without line search can exhibit residual growth at several first iterations and then still converge quadratically.
\par
     The main problem of the Newton method is the need for an initial estimate sufficiently close to the solution. In some rather simple cases, the method can converge even with arbitrary constant hydraulic head distribution given as the initial estimate. However, for more difficult problems including heterogeneous domains with complex geometry and contrasting and anisotropic media, the task of choosing a good initial estimate is not trivial.

Since the Newton method for discretized equation \eqref{eq:RichStat} generally doesn't converge even with line search, another strategies are desired. Here two of them are described: the pseudo-transient method \cite{farthing2003efficient}\cite{anuprienko2018modeling}, and a continuation method based on varying nonlinearity of the equation \eqref{eq:RichStat}. Application of the pseudo-transient method to the steady-state Richards equation and its comparison with other solvers can be also found in \cite{farthing2003efficient}.

\subsection{Pseudo-transient method}
This method is based on solution of transient equation \eqref{eq:Richards} until the solution reaches steady state. More details on numerical solution of transient problem in GeRa can be found in \cite{anuprienko2018modeling}, including comparison of Picard and Newton methods and brief mention of the pseudo-transient method. 
\par
Discretization in time is done using fully implicit scheme, which at each time step results in a system of nonlinear equations, that is again solved with the Newton method. Initial guess is the solution from the previous time step, which means better convergence, since with sufficiently small time step size that initial guess becomes good enough. Therefore, in transient case convergence of the Newton method usually can be achieved without sophisticated techniques such as line search. Similar results were reported by other researchers \cite{farthing2003efficient}.
\par
After solution of the nonlinear system at a time step the residual of the system \eqref{eq:NonlinSys} is evaluated and convergence check for the steady-state problem is performed. Then the time step size $\Delta t$, which may vary in user-specified limits, is increased to $\Delta t = 1.5\Delta t$ if the number of nonlinear iterations at the time step did not exceed user-specified number $numit_{inc}$, and in case of failure of the Newton method the time step is decreased to $\Delta t = 0.5\Delta t$. 
\par
      The pseudo-transient method has the following drawbacks:
\begin{itemize}
\item the time which the system takes to achieve steady-state can be very large and is not known, although it can be estimated by calculation on a coarse grid;
\item severe time step restrictions can be experienced in case of highly heterogenous and relatively dry soils;
\item the method can be sensitive to initial hydraulic head distrubution.
\end{itemize}

\subsection{Nonlinearity continuation method}
	The method proposed in this paper is based on continuation methodology. Continuation methods, when applied for nonlinear systems obtained by discretization partial differential equation in form \eqref{eq:NonlinSys}, solve the systems by incrementally approaching the solution through series of intermediate problems. Methods of this type have been successfully applied for numerical solution of partial differential equations in computational fluid dynamics \cite{continCFD} and nonlinear solid mechanics \cite{continMechNewt}\cite{continMechAdap}.
\par
      As described above, the main problem when trying to solve discretized equation \eqref{eq:RichStat} with the Newton method is to provide a good initial guess. The idea of our nonlinearity continuation method is to find such initial guess as a solution of some "simpler" problem or a sequence of such problems. 
\par 
A \textit{continuation parameter} $q$ is introduced which controls "degree of nonlinearity" of the problem. Namely, an equation is introduced:
\begin{equation}\label{eq:ContinEq}
-\nabla \cdot (\mathcal{K}(h, q)\mathbb{K}\nabla h) = Q,
\end{equation}

where $\mathcal{K}(h, q)$ is some function such that 
\begin{equation}\label{eq:ContinFuncProp}
\mathcal{K}(h, 1) \equiv K_r(h), ~~~~ \mathcal{K}(h, 0) \equiv 1.
\end{equation}
\par
There are several possibilities to choose such function $\mathcal{K}(h, q)$ of which two are tested:

\begin{equation}\label{eq:ContinPow}
\mathcal{K}_{pow}(h, q) = (K_r(h))^q
\end{equation}
and
\begin{equation}\label{eq:ContinLin}
\mathcal{K}_{lin}(h, q) = 1 + q \cdot (K_r(h) - 1).
\end{equation}

\par
     Combining properties \eqref{eq:ContinFuncProp} with equation \eqref{eq:ContinEq} one can check that with $q = 1$ equation \eqref{eq:ContinEq} is the original equation \eqref{eq:RichStat} and with $q = 0$ it becomes a simple linear equation, to solve which only a linear system has to be solved. Now, with the solution of the linear problem obtained, an attempt may be made to solve "more nonlinear" problem with $q > 0$ with the Newton method and solution of the linear problem taken as initial guess. The algorithm of incremental increasing of the continuation parameter $q$ to $1$ is described in algorithm \ref{alg:Contin}.

\begin{algorithm}[H]\label{alg:Contin}
\SetAlgoLined
    $q = 0$\;
    Solve linear problem $-\nabla \cdot (\mathbb{K}\nabla h) = Q$\;
    $\Delta q_{last} = 1$\;
    
    \While{$q < 1$}{
    	$\Delta q = \min(1-q~,2\cdot\Delta q_{last})$\;
    	\While{$\Delta q > 10^{-4}$}{
            Try to solve $-\nabla \cdot (\mathcal{K}(h, q+\Delta q)\mathbb{K}\nabla h) = Q$ with Newton\;
            \eIf{solved}{
                save $h$\;
                $\Delta q_{last} = \Delta q$\;
                break\;
            }{
                $\Delta q  = \Delta q /2$\;
            }
        }
        \If{didn't find $\Delta q$}{
            failed\;
            stop\;
        }
    }
 \caption{Nonlinearity continuation method}
\end{algorithm}

\subsection{Implementation details}
GeRa is written in C++ with its MPI-parallel computational part based on INMOST \cite{danilov2016parallel}, a parallel computation platform that provides tools for working with meshes consisting of arbitrary polyhedral cells, mesh partitioners and linear solvers as well as some other useful features. Automatic differentiation capabilities of INMOST are used to construct the Jacobian matrix. For solution of system \eqref{eq:JacSys} we use algorithms based on Bi-CGSTAB \cite{van1992bi}:
\begin{itemize}
	\item ILU(3)-preconditioned solver from PETSc package \cite{balay2019petsc};
	\item Inner\_MPTILUC, an internal INMOST solver based on second order Crout-ILU with inverse-based condition estimation and maximum product transversal reordering as preconditioner.
\end{itemize}
\par Of these solvers the latter is more robust and is used when PETSc solver fails to converge, such behavior was noted for some complex problems, especially when using MPFA scheme.

\section{Numerical experiments}
This section includes a number of tests comparing approaches described above on model and real-life problems.

\subsection{Problems setup}
Numerical experiments were performed for 3 different problems, of which the first is a two-dimensional model problem describing groundwater flow through a dam and the remaining two problems are real-life problems describing groundwater flow in sites that will be further denoted as site A and site B.  
\par
The first test, a model problem taken from \cite{polubarinova2015theory}, describes two-dimensional groundwater flow through a dam. The dam is a 10 m $\times$ 10 m square composed of homogeneous material with isotropic hydraulic conductivity $K = 0.864$ m/day. The left boundary has prescribed hydraulic head value $h = 10$ m, while the right boundary has prescribed hydraulic head value $h = 2$ m up to the point where $z = 2$ m, above that point a seepage boundary condition is imposed. Top and bottom boundaries have zero-flux boundary conditions. Some aspects of numerical solution of the dam problem are described in \cite{anuprienko2018modeling}. The problem was considered as a quasi-two-dimensional: a third dimension was added, along which only one layer of cell was present, and zero-flux boundary conditions were imposed on boundaries orthogonal to the added dimension axis.
\par
Real-life site A is a domain of approximately 64 km$^2$ with 3 geological layers, 7 different media and presence of lakes and rivers. All media have diagonal anisotropic hydraulic conductivity tensor of the form $\mathbb{K} = \text{diag}\{K, K, 0.1K\}$ with K varying from 0.011 to 4.76 m/day. This problem was also mentioned in \cite{anuprienko2018modeling}.
\par
Real-life site B is a domain of approximately 100 km$^2$ with 4 geological layers, 4 different media and presence of lakes and rivers. All media have diagonal anisotropic hydraulic conductivity tensor of the form $\mathbb{K} = \text{diag}\{K, K, 0.1K\}$ with K varying from 0.14 to 16 m/day.
\subsection{Comparison of pseudo-transient and nonlinearity continuation methods on the dam problem}
Computations were carried out on two cubic grids with 1600 cells (cell size 0.25 m) and 10000 cells (cell size 0.1 m). Discretization scheme was TPFA with upwind approximation of $K_r$. The following parameters were used: $\varepsilon_{rel} = 10^{-5}$, $\varepsilon_{abs} = 10^{-5}$, $maxit = 25$ and (in pseudo-transient method) $numit_{inc} = 15$.
\par
\begin{figure}
	\centering
	\includegraphics[width=0.4\textwidth]{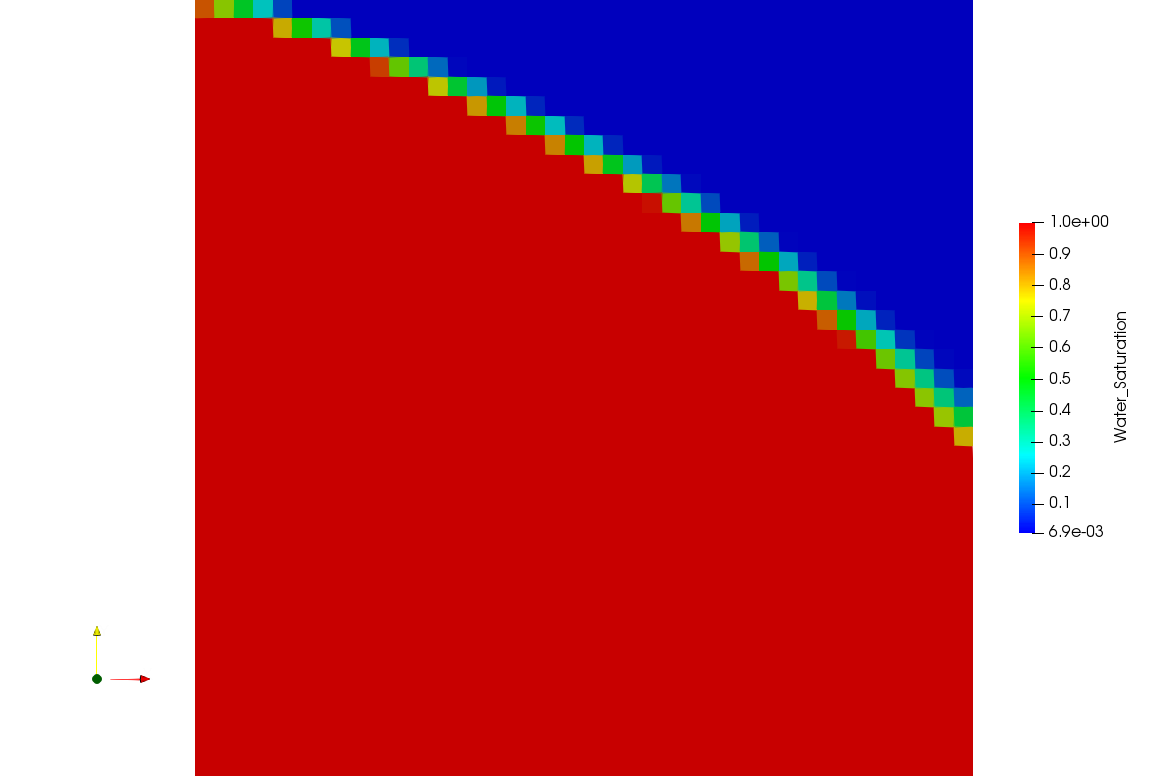}
	\includegraphics[width=0.4\textwidth]{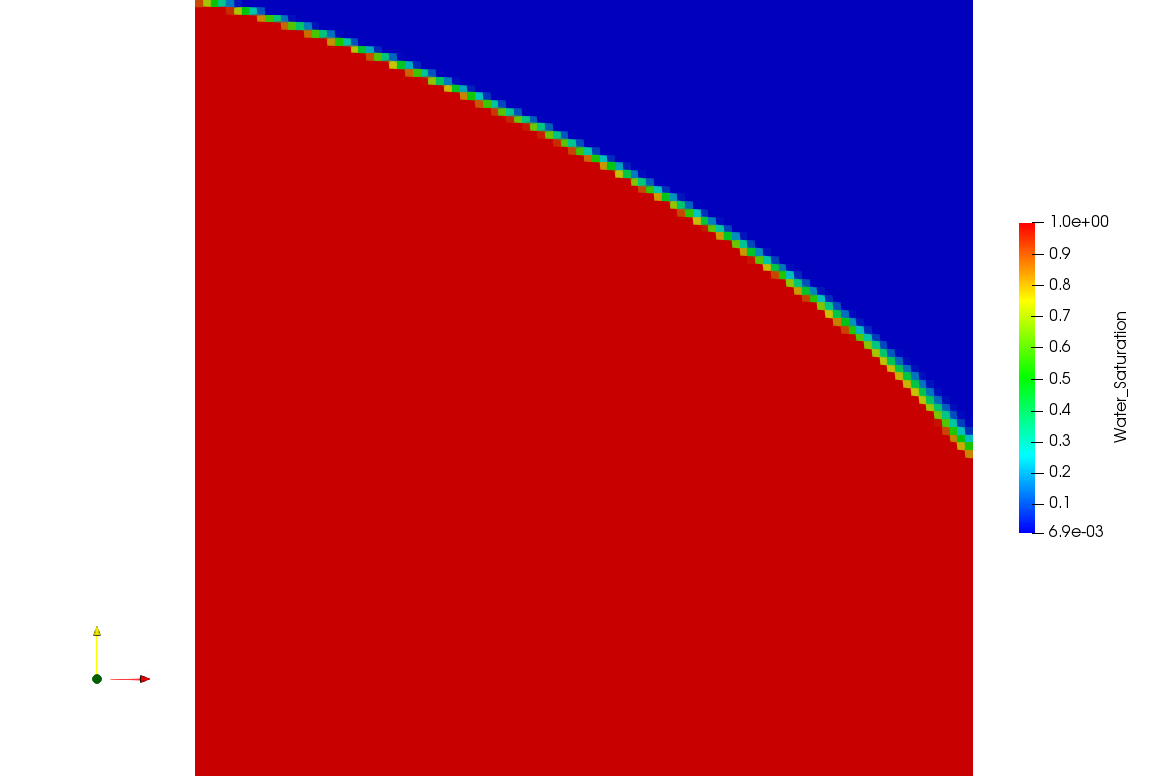}
	\caption{Saturation profiles for the dam test on cubic grids}\label{pic:damTPFA}
\end{figure}

Computation time, number of steps (time steps for pseudo-transient method and continuation steps for continuation method) and total number of Newton iterations are presented in tables \ref{tab:dam1600} and \ref{tab:dam10000}; saturation profiles are presented in figure \ref{pic:damTPFA}. Experiments showed that on both meshes continuation method took only 1 step to converge and much smaller number of Newton iterations compared to pseudo-transient method, leading to significantly faster calculations.

\begin{table}
	\begin{center}
		\begin{tabular}{ c| c| c| c }
			& $T_{comp}$, s & \# of steps & \# of Newton iterations \\\hline 
			Pseudo-transient & 2.65          &  8          &  49   \\  
			Continuation     & 0.85          &  1          &   8      
		\end{tabular}
		\caption{Comparison of pseudo-transient and continuation methods for the dam test, mesh with 1600 cells, TPFA}
		\label{tab:dam1600}
	\end{center}
\end{table}

\begin{table}
	\begin{center}
		\begin{tabular}{ c| c| c| c }
			& $T_{comp}$, s & \# of steps & \# of Newton iterations \\\hline 
			Pseudo-transient & 18.9          &  8          &  63   \\  
			Continuation     &  6.2          &  1          &  13      
		\end{tabular}
		\caption{Comparison of pseudo-transient and continuation methods for the dam test, mesh with 10000 cells, TPFA}
		\label{tab:dam10000}
	\end{center}
\end{table}

\subsection{Comparison of pseudo-transient and nonlinearity continuation methods on the site A}
Computations were carried out on triangular prismatic mesh of 171450 cells. The following parameters were used: $\varepsilon_{rel} = 10^{-5}$, $\varepsilon_{abs} = 10^{-6}$, $maxit = 25$ and (in pseudo-transient method) $numit_{inc} = 15$.
\par
To compare pseudo-transient and continuation methods the problem was solved using TPFA scheme; additional calculation was performed with continuation method and MPFA scheme to show that with continuation method it is possible to solve problem with more accurate and complex MPFA scheme still faster than with pseudo-transient method and simple TPFA. All calculations employed central approximation of $K_r$.
\par
Computation time, number of steps (time steps for pseudo-transient method and continuation steps for continuation method) and total number of Newton iterations are presented in table \ref{tab:obj1_171k}.

\begin{table}
	\begin{center}
		\begin{tabular}{ c| c| c| c }
		                         	& $T_{comp}$, s& \# of successful (failed) steps & \# of Newton iterations \\\hline 
			Pseudo-transient, TPFA & 868.4         &  65(10)     & 377   \\  
			Continuation, TPFA     & 124.4         &   1(0)      &  10   \\  
			Continuation, MPFA     & 444.3         &   1(0)      &  13           
		\end{tabular}
		\caption{Comparison of pseudo-transient and continuation methods for site A, mesh with 171450 cells}
		\label{tab:obj1_171k}
	\end{center}
\end{table}


\subsection{Effect of line search and comparison of linear and power functions in continuation method}
Line search in Newton method is a useful tool which often allows to increase continuation parameter $q$ from 0 to 1 in a single step, as can be seen in tables \ref{tab:dam1600}, \ref{tab:dam10000} and \ref{tab:obj1_171k}. In these cases, there is no difference between using power \eqref{eq:ContinPow} and linear \eqref{eq:ContinLin} functions as continuation functions. To show the difference, either line search should be omitted or more complex problems should be regarded.
\par
\begin{figure}
	\centering
	\includegraphics[width=0.6\textwidth]{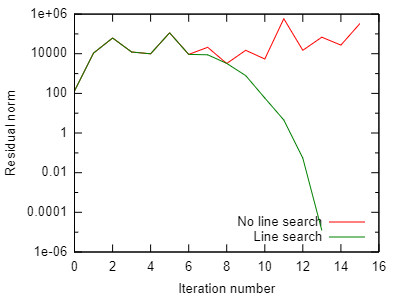}
	\caption{Effect of line search for site A. A single iteration with line search ($\Omega = 0.25$ at iteration 7) allows Newton method then to converge quadratically, otherwise Newton doesn't converge}\label{pic:ls_effect}
\end{figure}

Computations for site A on the same mesh were performed with no line search in the Newton method and maximal iterations number decreased to $maxit = 15$. Computation time, number of steps and total number of Newton iterations are presented in table \ref{tab:obj1_linpow}. Effect of line search on Newton performance in contiuation method is shown in figure \ref{pic:ls_effect}; a single iteration needed application of line search, but that drastically changed overall convergence. When using TPFA, there is little difference between computational effort with linear and power functions; however, when using MPFA, linear function leads to one more failed continuation step and subsequently slows down the computation. 
\par
\begin{figure}
	\centering
	\includegraphics[width=0.8\textwidth]{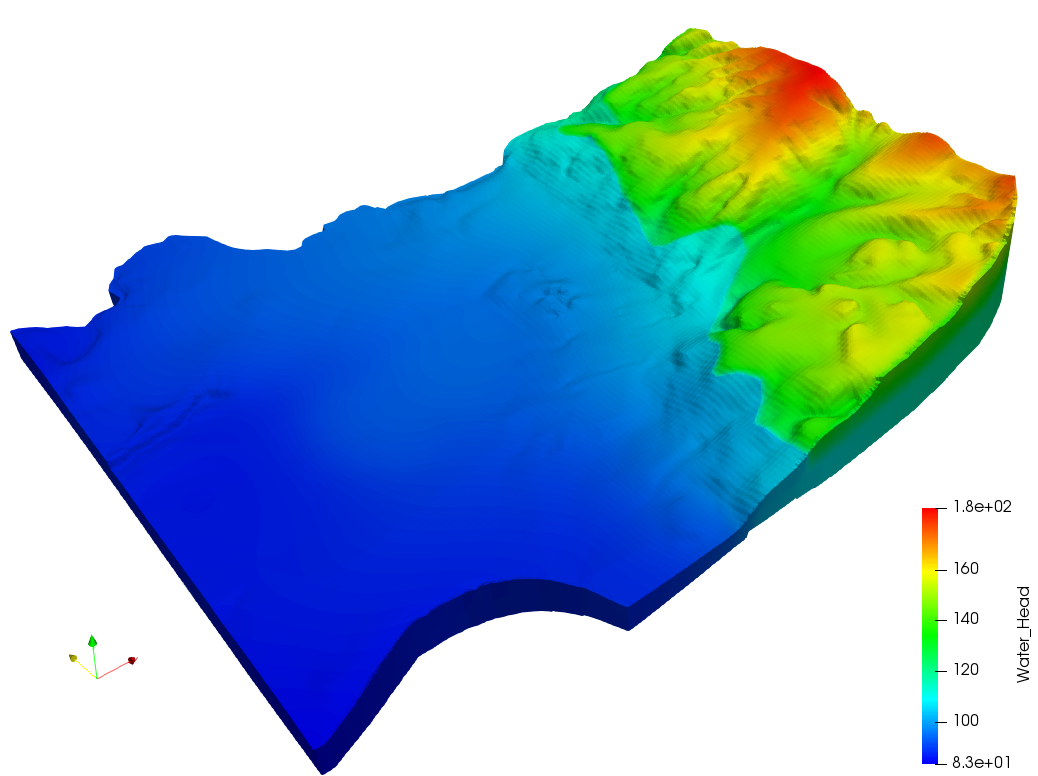}
	\includegraphics[width=0.8\textwidth]{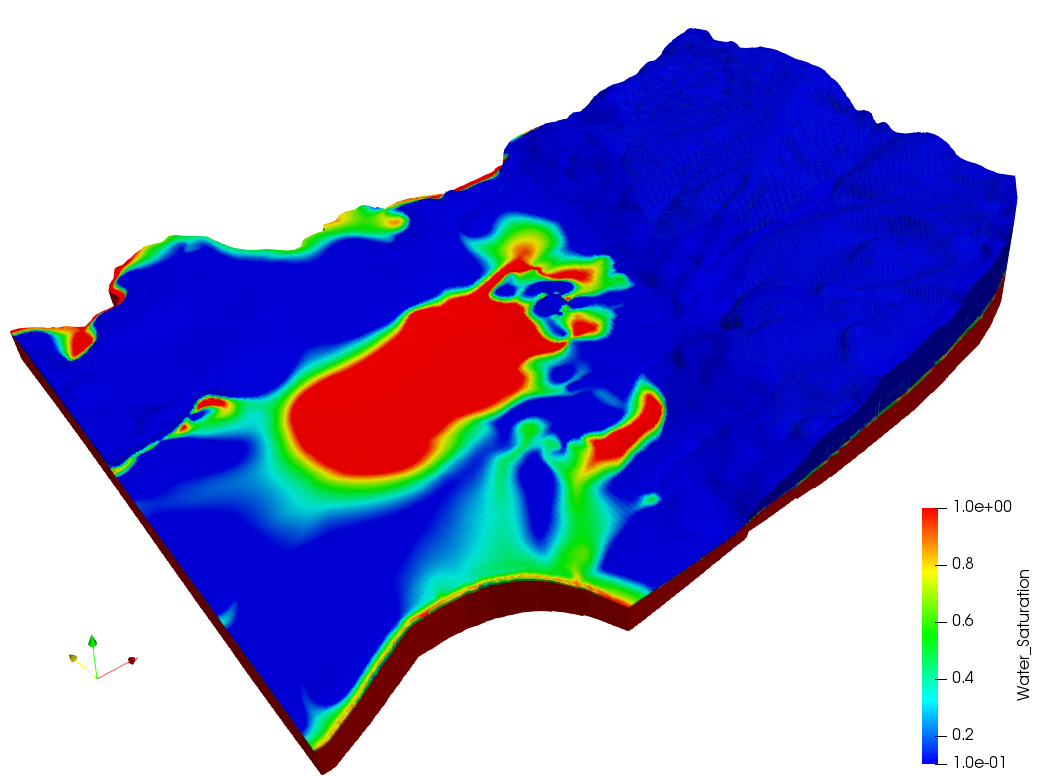}
	\caption{Hydraulic head and saturation, site A, mesh with 12 million cells, MPFA}\label{pic:12m_head_sat}
\end{figure}
\begin{figure}
	\centering
	\includegraphics[width=0.45\textwidth]{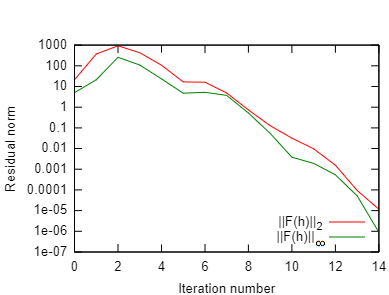}
	\includegraphics[width=0.45\textwidth]{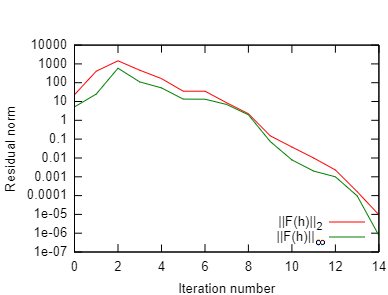}
	\caption{Example of convergence behavior of the Newton method with line search. On the left is linear continuation function, increasing $q$ from 0.75 to 1, on the right is power continuation function, increasing $q$ from 0.5 to 1. Site A, mesh with 12 million cells, MPFA}\label{pic:conv}
\end{figure}

To further investigate difference between linear and power functions, another tests were performed on a fine mesh with 12 million cells and MPFA discretization. The computations were carried out using 140 cores of the INM RAS cluster \cite{cluster}. The results are shown in table \ref{tab:obj1_cluster}, the hydraulic head and saturation distributions are presented in the figure \ref{pic:12m_head_sat}. Example of convergence behavior of the Newton method with line search  is shown in figure \ref{pic:conv}. In case of such fine mesh resolution nonlinearity continuation required more than one step even with line search in the Newton method. Unlike in the case of coarse grid, the power function in the continuation method lead to smaller number of steps in the continuation method and smaller total number of Newton iterations, but solution time was comparable to that of the linear function. 

\begin{table}
	\begin{center}
		\begin{tabular}{ c| c| c| c }
			& $T_{comp}$, s & \# of successful (failed) steps & \# of Newton iterations \\\hline 
			Power, TPFA    & 325.8         &   2(1)      &  32   \\  
			Power, MPFA    & 1069.3        &   2(1)      &  32   \\  
			Linear, TPFA   & 347.5         &   2(1)      &  34   \\
			Linear, MPFA   & 1445.1        &   3(2)      &  51   
		\end{tabular}
		\caption{Comparison of power and linear continuation functions, site A, mesh with 171450 cells, Newton method without line search}
		\label{tab:obj1_linpow}
	\end{center}
\end{table}

\begin{table}
	\begin{center}
		\begin{tabular}{ c| c| c| c }
			& $T_{comp}$, s & \# of successful (failed) steps & \# of Newton iterations \\\hline 
			Linear  &  8776         &  3(2)       &  55   \\
			Power   &  8374         &  2(1)       &  30   
		\end{tabular}
		\caption{Results of parallel computations for site A, mesh with 12 million cells, MPFA}
		\label{tab:obj1_cluster}
	\end{center}
\end{table}

\par
Site B turned out to be harder problem in the sense that even with line search in Newton method continuation method takes several steps even on relatively coarse grids. The following parameters: $\varepsilon_{rel} = 10^{-6}$, $\varepsilon_{abs} = 10^{-6}$, $maxit = 25$, central approximation of $K_r$ and a trinagular prismatic mesh of 56820 cells were used. Computation time, number of steps and total number of Newton iterations are presented in table \ref{tab:obj2_linpow}. When using TPFA, there is again little difference in computational time with linear and power functions, although power functions takes twice more steps; when using MPFA, linear function performs better leading to less continuation steps and Newton iterations.

\begin{table}
	\begin{center}
		\begin{tabular}{ c| c| c| c }
			               & $T_{comp}$, s & \# of successful (failed) steps & \# of Newton iterations \\\hline 
			Power, TPFA    &  171.1        & 10(9)       & 133   \\  
			Power, MPFA    &  732.3        &  7(8)       & 117   \\  
			Linear, TPFA   &  163.9        &  5(4)       & 104   \\
			Linear, MPFA   &  393.6        &  5(4)       &  61   
		\end{tabular}
		\caption{Comparison of power and linear continuation functions, site B, mesh with 56820 cells, Newton method with line search}
		\label{tab:obj2_linpow}
	\end{center}
\end{table}

\begin{figure}
	\centering
	\includegraphics[width=0.8\textwidth]{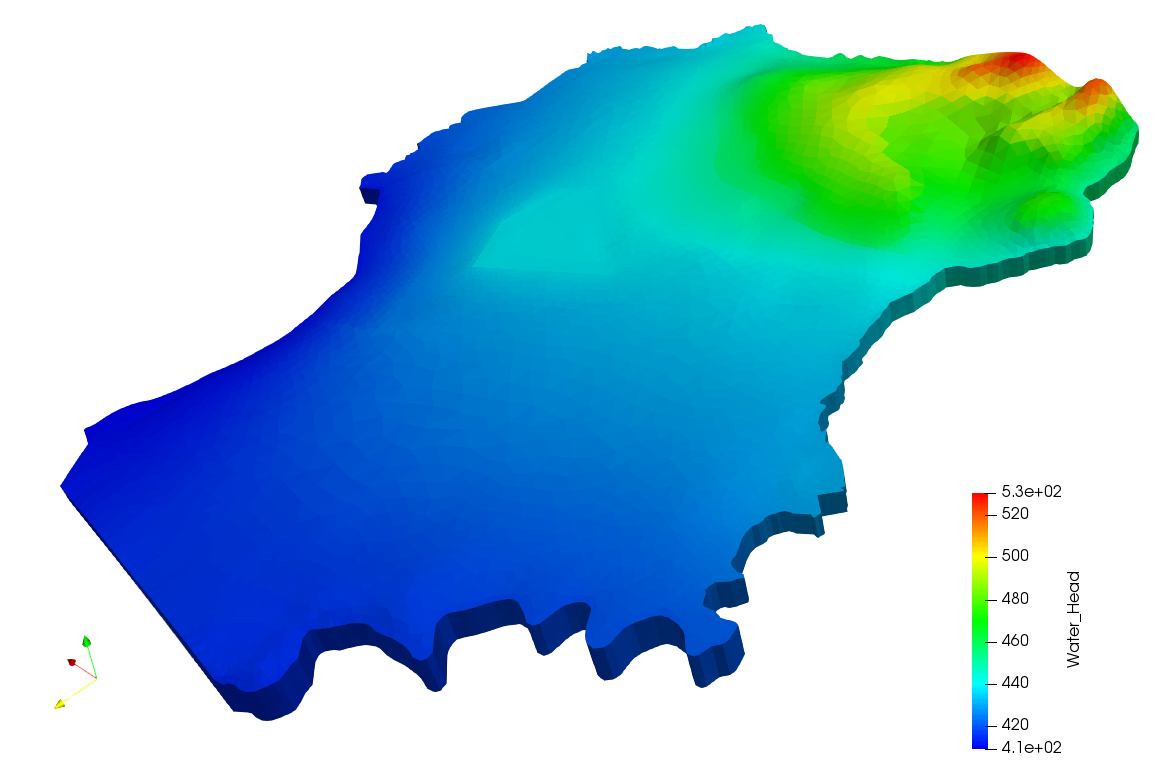}
	\includegraphics[width=0.8\textwidth]{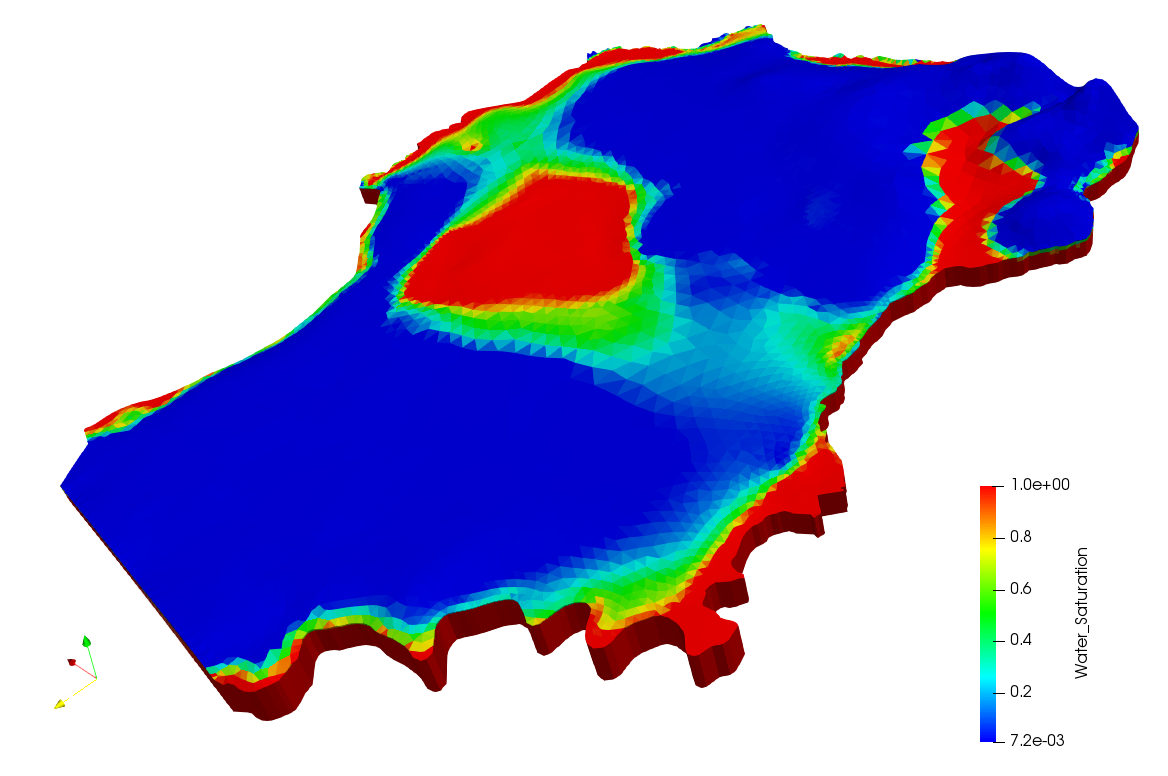}
	\caption{Hydraulic head and saturation distributions for site B, mesh with 56820 cells, MPFA}\label{pic:obj2_58k}
\end{figure}

\section{Conclusion}
A nonlinearity continuation method for solution of nonlinear systems arising from discretization of steady-state Richards equation for modeling  groundwater flow in variably saturated conditions is proposed. Through solution of a series of problems with increasing nonlinearity the method gets a suitable initial estimate required by the Newton method. The continuation procedure is implemented in the GeRa software package.
\par Combined with line search for improved convergence of Newton method, the continuation method allows for solution of complex problems with MPFA discretization scheme consuming less time than previously used pseudo-transient method does with simpler TPFA discretization scheme.
\par Two continuation functions, linear and power, were compared on two real-life problems; however, the tests didn't show that any of the two is superior than the other, with their performance in the continuation method depending on the problem, discretization scheme and mesh size.

\section*{Acknowledgements}
This work was partially (namely development and testing of line search technique in Newton method) supported by Russian Science Foundation through the grant 18-71-10111 and partially (namely preparation and execution of parallel computations) supported by the Russian Academy of Sciences Research program No. 26 "Basics of algorithms and software for high performance computing".

\bibliographystyle{elsarticle-num} 
\bibliography{jcam2020.bib}





\end{document}